\documentclass[10pt]{article}

\usepackage[utf8]{inputenc}
\usepackage{amscd,amssymb}
\usepackage{amsthm,amsmath}
\usepackage{graphicx}
\usepackage{enumerate,amssymb,setspace}
\usepackage{amsfonts,amsthm,amscd}
\usepackage{mathtools}
\usepackage{cite}
\usepackage{mathtools}
\usepackage[all]{xy}
\usepackage{hyperref}

\usepackage{caption}
\usepackage[]{subcaption}

\DeclareCaptionLabelFormat{opening}{}

\newtheorem{theorem}{Theorem}[section]
\newtheorem{lemma}[theorem]{Lemma}
\newtheorem{prop}[theorem]{Proposition}
\newtheorem{corollary}[theorem]{Corollary}
\newtheorem{conjecture}[theorem]{{Conjecture}}

\newtheorem{problem}[theorem]{{Problem}}

\newtheorem{definition}[theorem]{{Definition}}
\newtheorem{claim}[theorem]{{Claim}}

\def\bclaim{\begin{claim}}
	\def\eclaim{\end{claim}}
\def\bdefin{\begin{definition}}
	\def\edefin{\end{definition}}
\def\bcor{\begin{corollary}}
	\def\ecor{\end{corollary}}
\def\bthm{\begin{theorem}}
	\def\ethm{\end{theorem}}
\def\bconj{\begin{conjecture}}
	\def\econj{\end{conjecture}}
\def\blem{\begin{lemma}}
	\def\elem{\end{lemma}}
\def\blemma{\begin{lemma}}
	\def\elemma{\end{lemma}}
\def\bprop{\begin{prop}}
	\def\eprop{\end{prop}}
\def\bremark{\begin{remark}}
	\def\eremark{\end{remark}}
\def\bprob{\begin{problem}}
	\def\eprob{\end{problem}}

\theoremstyle{remark}
\newtheorem{remark}[theorem]{Remark}

\newcommand{\Pic}{\mathrm{Pic}}

\newcommand{\rk}{\mathrm{rk}}

\newcommand{\Supp}{\mathrm{supp}}
\newcommand{\mult}{\mathrm{mult}}

\makeatletter\@addtoreset{equation}{section} \makeatother

\font\sml=cmr6  
\def\al{\alpha}
\def\be{\beta}

 \def\eps{\epsilon}
\def\K{K\"ahler } 
\def\h#1{\hbox{#1}}

\def\q{\quad} 
\def\ra{\rightarrow}

\def\Fn{\mathbb{F}_n}
\def\Pic{\operatorname{Pic}}
\def\sm{\setminus}

\newcommand{\PP}{{\mathbb P}} \newcommand{\RR}{\mathbb{R}}
\newcommand{\QQ}{\mathbb{Q}} 
 \newcommand{\NN}{{\mathbb N}}
\newcommand{\FF}{{\mathbb F}}
\def\Fn{\FF_n}

\def\beq{\begin{equation}}
\def\eeq{\end{equation}}
\def\bpf{\begin{proof}}
	\def\epf{\end{proof}}
\def\bremark{\begin{remark}}
	\def\eremark{\end{remark}}

\def\eaeq{\end{aligned}}
\def\baeq{\begin{aligned}}
\def\saldp{strongly asymptotically log del Pezzo }
\def\Saldp{Strongly asymptotically log del Pezzo }
\def\saldps{strongly asymptotically log del Pezzos }

\def\aldp{asymptotically log del Pezzo }

\def\aldps{asymptotically log del Pezzos }

\def\saldpno{strongly asymptotically log del Pezzo}

\def\aldpno{asymptotically log del Pezzo}

\def\alf{asymptotically log  Fano }

\def\ldps{log del Pezzos }
\def\ldpno{log del Pezzo}

\def\lb{\label}
\def\er{\eqref}
\def\noi{\noindent}

\newcommand\blfootnote[1]{%
	\begingroup
	\renewcommand\thefootnote{}\footnote{#1}%
	\addtocounter{footnote}{-1}%
	\endgroup
}

\title{
On
the body of ample angles
of asymptotically log Fano varieties}

\author{Paolo Cascini, Jesus Martinez-Garcia,
Yanir A. Rubinstein}
\begin{document}
	\bibliographystyle{amsalpha}
	\maketitle
	\begin{abstract}
		In dimension two, we reduce the
		classification problem for asymptotically log Fano pairs to 
		the problem of determining generality conditions on certain blow-ups.
		In any dimension, we prove the rationality of the body of ample angles 
		of an asymptotically log Fano pair, i.e., these convex bodies are 
		always rational polytopes.
		
	\end{abstract}

	\blfootnote{This note grew from  interactions between J.M.G. and Y.A.R, Y.A.R.'s lecture, and the ensuing conversations between P.C. and Y.A.R. at the conference on
		``Birational Geometry, K\"ahler--Einstein Metrics and Degenerations" that took
		place in November 2019. Thanks go to J. Park and POSTECH for the 
		excellent conference and hospitality.
		Thanks to I.A. 
		Cheltsov for co-organizing the conference as well as many helpful discussions.	
		The research of P.C. was supported by an EPSRC fellowship. 
		The research of Y.A.R. was supported by
		NSF grants DMS-1515703,1906370 and the Rosi \& Max Varon Visiting
		Professorship (Fall 2019 and Spring 2020) at the Weizmann Institute of Science to which he is
		grateful for the excellent research conditions.
		
		MSC subject codes: 14J45, 14J26 (primary), 14J10, 14E05 (secondary).
		
	}
	
	\section{Introduction}
	Asymptotically log Fano pairs, introduced by Cheltsov--Rubinstein \cite{CR}, generalizing work of \cite{Maeda}, have received attention in the last decade within the theory of K-stability and the Calabi problem \cite{CR2, CRZ, Fuj1, Fuj3}. We believe asymptotically log Fanos to be interesting objects in their own right. Indeed, while they belong to an infinite number of deformation families, one can expect a classification to be achievable. Alas, so far they have only been systematically studied under two special assumptions: the boundary having only one component or its generalization: the \emph{strong} regime. 
	
	Our goal in this note is to expand the knowledge on the birational geometry of asymptotically log Fano pairs by considering two separate problems.
	First, we discuss
	some aspects of the classification
	of asymptotically log Fano pairs in dimension two, and
	reduce such a classification to the understanding
	of generality conditions on the location of points, possibly
	infinitely near, that are blown-up on a curve in  a rational surface
	with Picard group of rank at most two.
	Second, we point out the rationality of certain
	convex bodies that arises in the study of asymptotically log Fano pairs in any dimension,
	making contact with recent advances on log-geography
	and the minimal model program.
	
	\subsection{The classification problem}
	
	A pair $(X,D)$ consisting of a smooth
	complex projective variety $X$ and a simple-normal-crossing (holomorphic) effective non-zero divisor 
	$D=\sum_{i=1}^rD_i$ in $X$ (where the $D_i$ are distinct irreducible hypersurfaces) 
	is called {\it asymptotically log
		Fano} in the sense of Cheltsov--Rubinstein \cite[Definition 1.1]{CR} if there exists a sequence
	$$\be(j)=(\be_1(j),\ldots,\be_r(j))
	\in(0,1]^r\cap\QQ^r, \,j\in\NN,$$ 
	converging to the
	origin such that 
	\beq
	\lb{aldpdefeq}
	-K_X-\sum_{i=1}^r(1-\be_i(j))D_i
	\q \h{is a $\QQ$-ample divisor, for each $j\in\NN$}.
	\eeq
	
	\noi
	A pair $(X,D)$ is called {\it strongly asymptotically log
		Fano} if \er{aldpdefeq} even holds in a semi-open sub-cube, not just
	along a sequence, i.e.,
	\beq
	\lb{saldpdefeq}
	-K_X-\sum_{i=1}^r(1-\be_i)D_i
	\q \h{is a $\QQ$-ample divisor,
		for all $\be\in(0,\eps]^r\cap\QQ^r$ for some $\eps\in(0,1]$}.
	\eeq
	Finally, $(X,D)$ is called {\it log
		Fano} if \er{saldpdefeq} actually holds in a closed sub-cube,
	i.e., for $\be\in[0,\eps]^r\cap\QQ^r$
	(to avoid confusion, we remark that some
	authors use this terminology to refer to rather more general
	objects 
	\cite[Definition 2.7]{CG}
	than Maeda's Definition \ref{definition:log-del-Pezzo-pair} below).
	In dimension two, we will denote a pair by $(S,C)$ 
	and refer to \er{aldpdefeq}--\er{saldpdefeq} as \emph{(strongly) \aldpno}.
	
\begin{figure}
		\centering
		\includegraphics[width=0.3\textwidth]{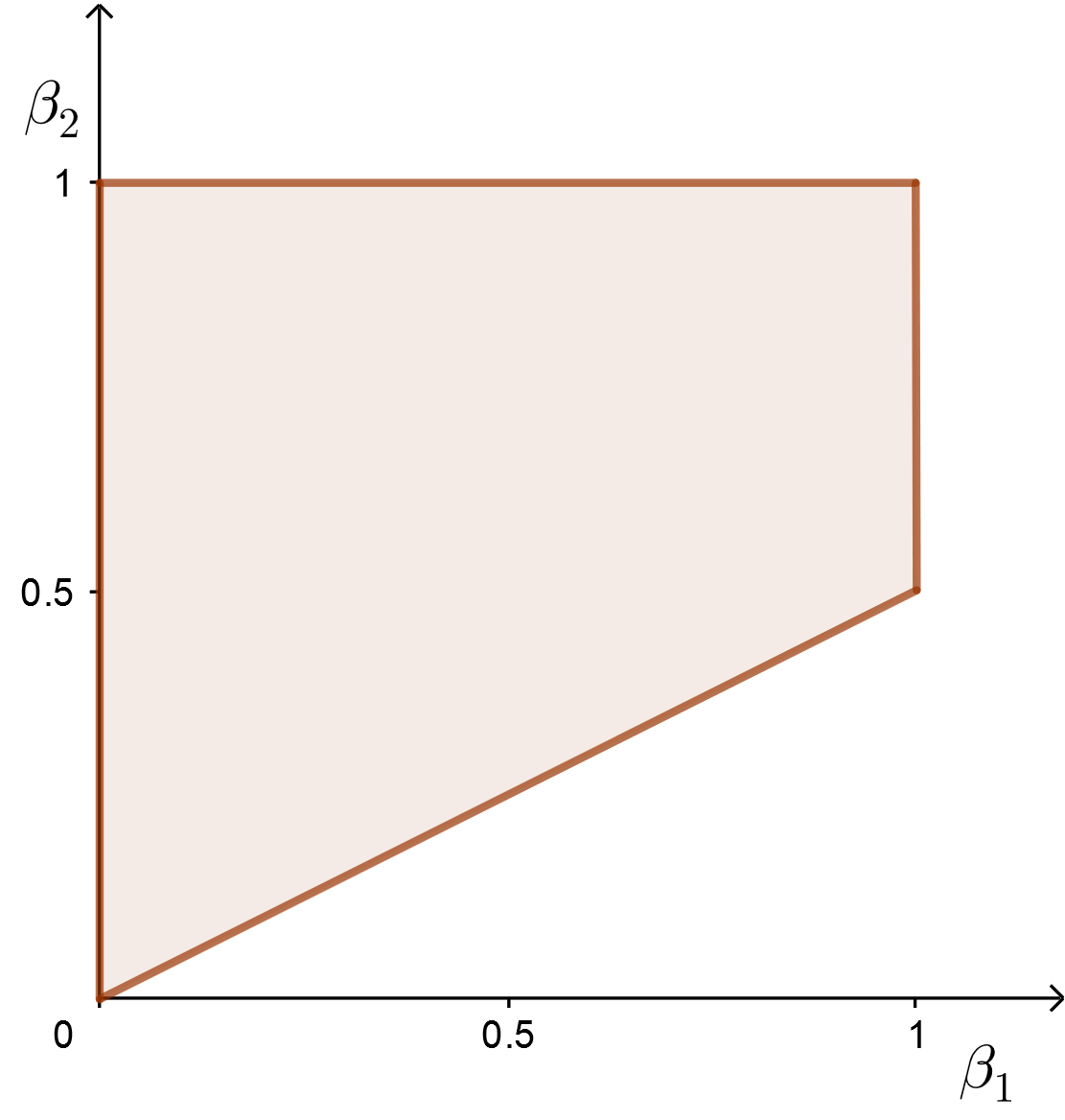}
		\caption{The body of ample angles for $(\FF_1,C_1+C_2)$ with
			$C_1=Z_1$ and $C_2$ smooth in $|-K_{\FF_1}-C_1|$ and intersecting
			$C_1$ transversally (see \S\ref{Hirbesubsec} for notation).}
		\label{Figure1}
	\end{figure}

	By definition, a log Fano pair $(X,D)$ is strongly asymptotically log
	Fano, and a strongly asymptotically log
	Fano pair $(X,D)$ is  asymptotically log
	Fano, but none of the reverse implications are true in general, already
	in dimension 2.
	Perhaps the simplest examples arise by considering anticanonical boundaries. If $S$ is del Pezzo and $C\sim-K_S$ is smooth
	then $(S,C)$ is \saldp but not \ldpno. The notions of \saldp and \aldp are actually equivalent
	when $C$ is smooth, 
	but as soon as $C$ contains two components they are not: let 
	$c\subset\PP^2$ be a cubic curve smooth away from a double point $p$ and let
	$S$ be the
	blow-up of $\PP^2$ at $p$. Let $C$ be the total transform of $c$, i.e.,
	$C=C_1+C_2$ with $C_1$ the inverse image of $p$, and
	with  $C_2$ the proper transform of $c$. In this case, \er{saldpdefeq}
	will not hold, in fact $-K_S-(1-\be_1)C_1-(1-\be_2)C_2$ will be ample if and only if $0<\be_1<2\be_2$ which defines a trapezoidal subset of the square
	$[0,1]^2$ (see Figure \ref{Figure1}).

	The following problem, raised by Cheltsov
	and one of us \cite{CR} (see \cite[\S8--9]{R14} for detailed 
	exposition), can be considered as a logarithmic 
	generalization of the folkore classification
	\cite{DelPezzo,Hitchin1975} of del Pezzo surfaces:
	
	\bprob
	\lb{mainprob}
	Classify \aldp surfaces.
	\eprob

	The classification of the subclass of {\it strongly} \aldps has
	been achieved \cite[Theorems 2.1, 3.1]{CR} (see also
	\cite{R20b}), however, as we will see below, the non-strongly
	regime is considerably harder to classify.
	Part of the motivation in \cite{CR}
	is the theory of
	canonical \K metrics with edge singularities
	\cite{JMR,MR} that is associated with such pairs.
	While Fano varieties exhibit remarkable
	finiteness properties in a given dimension,
	this is no longer the case for \alf varieties. In fact, 
	already 
	Maeda's classification of log del Pezzo surfaces, that we review below, exhibits
	infinitely many pairs.
	However,
	it is precisely the case of nonzero 
	$\beta_i$'s, interpreted geometrically
	by the existence of \K metrics with
	cone angle $2\pi\be_i(j)$ along $C_i$,
	that is of major interest in complex geometry, see \cite{R14} for
	a detailed survey and references. Moreover, as the recent works 
	\cite{CR,R14,CR2,CRZ,CMG, GMGS,Fuj1, Fuj3 ,R19,FujZhangetal,CZ} show, the notion of \alf varieties
	is interesting also purely from an algebraic geometry viewpoint. Finally, the family of \saldps is already
	much more vast and rich geometrically than the class of log del Pezzos, and, as we try to explain in this
	note, the class of \aldps is, in a sense, yet an order of magnitude more vast.
	Our first goal in this note will be to illustrate with concrete examples some of
	the difficulties in solving Problem \ref{mainprob} and 
	give some first steps in such a classification. More precisely, 
	we will show that all \aldp pairs arise from pairs with Picard rank $2$ in a very explicit way
	(Proposition \ref{partialclassifY})
	and classify the latter (Proposition \ref{prop:classification-rk2-anticanonical}). Classifying
	all \aldp pairs remains a difficult open problem.

	\subsection{The body of ample angles}
	
	One can collect all ``admissible" coefficients 
	in the sense of \er{aldpdefeq}
	(equivalently interpreted as angles for which there exists a
	\K edge metric with angles $2\pi\be_i$ along $C_i$ and
	with positive Ricci curvature on $S\sm C$ \cite{JMR,R14,MR,GP})
	into a single body previously introduced by one of us  \cite[Definition 3.1]{R19}:
	
	\begin{definition}
		The set
		\label{def:AAEq}
		\begin{equation}
		\begin{aligned}
		\label{eq:AAEq}
		\h{\rm AA}(X,D):=\Big\{\be=(\be_1,\ldots,\be_r)\in(0,1)^r\,:\, 
		\h{$-K_X-\sum_{i=1}^r(1-\be_i)D_i$ is ample} \Big\}
		\end{aligned}
		\end{equation}
		is called the \emph{body of ample angles} of $(X,D)$.
	\end{definition}

	The problem of determining whether a given 
	pair $(X,D=\sum_{i=1}^rD_i)$ is asymptotically log Fano
	amounts to determining whether
	$$
	0\in\overline{\h{\rm AA}(X,D)}.
	$$ 
	Thus, this set is a fundamental object in the
	study of \alf varieties.
	
	When $(X,D)$ is strongly \alf this body is simply a cube near the origin,
	in particular locally polyhedral near the origin.

It is easy to see that the body of ample angles of an \alf variety is a convex body \cite[Lemma 3.3]{R19}. A natural 
	question is therefore:
	
	\begin{problem}
		\label{polytopeProb}
		Let $(X,D)$ be asymptotically log Fano.
		Is $\overline{\h{\rm AA}(X,D)}$ polyhedral?
	\end{problem}
	
	Of course, Problems \ref{mainprob} and \ref{polytopeProb} are related:
	if one had a classification as in Problem \ref{mainprob}, 
	then it would likely be possible
	to glean from it explicitly the bodies of ample angles, which would
	resolve Problem \ref{polytopeProb} in dimension 2. In fact,
	as shown by one of us, it turns out that Problem \ref{polytopeProb}
	can be solved affirmatively in dimension 2 by describing rather explicitly the linear inequalities cutting out the body of ample angles, however this approach specifically hinges on the Nakai--Moishezon criterion and other two-dimensional facts and does not generalize to higher dimensions \cite{R20}.
	Here, we resolve Problem \ref{polytopeProb} in
	any dimension. In fact, we prove a stronger statement that does not assume $(X,D)$ is
	asymptotically log Fano (i.e., $0\in\overline{\h{\rm AA}(X,D)}$) but merely that 
	$\h{\rm AA}(X,D)$ is non-empty (note that by openness of ampleness the former implies
	the latter).

	\bthm
	\lb{CasciniRThm}
	$\overline{\h{\rm AA}(X,D)}$ is either empty or a rational polytope,
	i.e., cut out by finitely-many linear inequalities with rational coefficients in $\be_1,\ldots,\be_r$.
	\ethm
	
	This result perhaps motivates renaming the body of ample
	angles the `angletope' or `ample-angletope'. 
	
	Already in dimension two---at least from the point of view of the Nakai--Moishezon 
	criterion---polyhedrality is not
	obvious: the aforementioned criterion involves infinitely-many linear
	inequalities (one for the intersection of \eqref{saldpdefeq} with each irreducible curve of $S$) as well as one quadratic inequality in the $\be_i$'s (corresponding to taking the self-intersection of \eqref{saldpdefeq}). The latter
	condition would not seem to affect the body being polyhedral {\it near the origin}. In fact, it only comes to play in the cases when a birational model $(s,c)$ with Picard rank $2$ of $(S,C)$, obtained by smooth contractions of curves intersecting $C$, satisfies $(K_s+c)^2=0$. Nevertheless, it is remarkable
	that even in these cases it can basically be replaced by a finite collection of linear
	inequalities, and that overall all the linear inequalities boil down to
	finitely-many. In fact,
	the proof of Theorem \ref{CasciniRThm}
	uses the recent progress on the study of Shokurov's log geography,  due to one of us together with Birkar, Hacon and M\textsuperscript cKernan \cite{BCHM} (see also \cite{CL}) that was a crucial ingredient in the proof of the finite generation of the canonical ring. 
	
	\subsection{Organization}
	The remainder of this note is organized as follows.
	In \S\ref{ldpsec}  we recall the classical result of Maeda on the classification of \ldps 
	and give a self-contained proof.
	In \S\ref{saldpsec}  we give a brief overview of the classification of \saldps
	following Cheltsov and one of us \cite{CR,R20b}.
	In \S\ref{aldpsec}  we explain which aspects of the approach in op. cit. extend to the
	setting of \aldps and which ones do not. Our first main result is a description of all \aldps as a subset of a particular class of pairs that are obtained as certain blow-ups of \aldps with small Picard group
	(Proposition \ref{partialclassifY}).
	In \S\ref{aldppic2sec}  we give an explicit classification of \aldps with small Picard group (Proposition \ref{prop:classification-rk2-anticanonical}). Propositions 
	\ref{partialclassifY} and \ref{prop:classification-rk2-anticanonical} combine to give
	a vast class of \aldps that are not \saldpno. We end in \S\ref{polytopesec}
	with a proof of Theorem \ref{CasciniRThm}.

	\section{Log del Pezzo surfaces}
	\lb{ldpsec}
	While the classification of log del Pezzo surfaces is rather simple, it serves
	to illustrate the simplest setting
	and the origin of Problem \ref{mainprob}.
	

	The definition of log Fano manifolds goes back to work of Maeda \cite{Maeda}.
	\begin{definition}
		\label{definition:log-del-Pezzo-pair} 
		Let $X$ be a smooth variety and let $D$ be a simple normal
		crossing divisor in $X$.
		We say that the
		pair $(X,D=\sum D_i)$ is \emph{ log Fano} if
		$-K_X-D$ is ample.
	\end{definition}
	\noindent
	
	In dimension 2, these are also called
	log del Pezzo surfaces.
	The motivation for the adjective
	``logarithmic", according to Maeda, is from the work of
	Iitaka on the classification of open algebraic varieties
	where logarithmic differential forms are used
	to define invariants of the pair. The open variety associated
	to $(X,D)$ is the Zariski open set $X\sm D$.

	\subsection{Some facts on Hirzebruch
		surfaces}
	\lb{Hirbesubsec}
	
	The main difference between del Pezzo
	surfaces and log del Pezzo surfaces is the
	appearance of Hirzebruch surfaces, $\Fn$. Let us recall some
	basic facts and establish some notation for $\Fn$.

	For each $n\ge0$,
	denote by $\Fn$ 
	the unique rational ruled surface
	whose Picard group has rank two and
	contains a unique (if $n>0$) smooth rational curve of self-intersection $-n$.
	We denote this curve by $Z_n$, and by $F$ we denote
	the class of an irreducible smooth rational curve such that $F^2=0$ and $F.Z_n=1$. 
	If $n=0$ when we refer to $Z_0$ and $F$ we
	intend that each is a fiber of a different projection to $\PP^1$.
	Hirzebruch surfaces are ruled toric surfaces  
	and applying adjunction yields \cite[Chapter 5, \S2]{Har77}
	\beq
	\label{KFnEq}
	-K_{S}\sim 2Z_n+(n+2)F
	\eeq
	Recall that every smooth irreducible curve in $|Z_n+nF|$
	(a `zero section') intersects each fiber transversally at a single point
	and does not intersect the `infinity section' $Z_n$.
	Any curve $C$ on $\FF_n$ satisfies 
	\beq\label{Fn-curves}
	C\sim aZ_n+bF
	\eeq
	with $a,b\in\NN\cup\{0\}$. This, combined with the Nakai-Moishezon criterion implies
	\beq\label{ampleFn}
	\h{$C$ is ample if and
		only if $a>0$ and $b>na$,}
	\eeq
	and
	\beq\label{nefFn}
	\h{$C$ is nef if and
		only if $a\ge0$ and $b\ge na$,}
	\eeq
	and furthermore,
	\beq\label{irreducibleFn}
	\h{$C$ is an
		irreducible curve
		only if $C=Z_n$ or $b\ge na\ge0$,}
	\eeq
	and under such conditions the class \eqref{Fn-curves} always contains an irreducible curve.

	\subsection{Maeda's classification}
	
	Maeda classified log del Pezzo pairs \cite[\S3.4]{Maeda}. Let us 
	review his classification as its rigidness serves as a contrast
	to the flexibility we demonstrate for \aldp pairs in Propositions \ref{partialclassifY} 
	and \ref{prop:classification-rk2-anticanonical} below.

	\begin{prop}\label{MaedaThm}
		Log del Pezzo surfaces $(S,C=\sum_{i=1}^rC_i)$ are classified as follows:
		\hfill\break
		(i) $S=\mathbb{P}^2$, $C_1$ is a line,%
		\hfill\break
		(ii)
		$S=\mathbb{P}^2$, $C_1,C_2$, are lines,%
		\hfill\break
		(iii)
		$S=\mathbb{P}^2$, and $C_1$ is a conic,%
		\hfill\break
		(iv)
		$S=\mathbb{F}_{n},\,n\in\NN\cup\{0\}$, $C_1=Z_n$,%
		\hfill\break
		(v)
		$S=\mathbb{F}_{n},\,n\in\NN\cup\{0\}$, $C_1=Z_n, C_2\in |F|$,
		\hfill\break
		(vi)
		$S=\mathbb{F}_1$, $C_1\in |Z_1+F|$,
		\hfill\break
		(vii)
		$S=\mathbb{P}^1\times\mathbb{P}^1$, $C_1$ is a $(1,1)$-curve.	
	\end{prop}
	
	\bpf
	By the Kodaira Vanishing Theorem and Kodaira--Serre duality, it follows that $S$ is rational (cf. \cite[\S3]{CR}). 
	When $S=\PP^2$,
	and since $-K_S\sim 3H$ we see the possibilities are 
	(i), (ii), or (iii).
	Assume now that $S=\Fn$. Denote
	$$
	C_i\in|a_iZ_n+b_iF|.
	$$
	Since $-K_S-C$ is ample and using \er{KFnEq} and \er{ampleFn}, we see that 
	$$
	\sum_i a_i\in\{0,1\}, \q 
	\sum_i b_i\in\{0,\ldots,n+1\}. 
	$$
	and that at most one component of $C$ is a fiber.
	Note that by \er{irreducibleFn} either $b_i\ge na_i>0$ or else $(a_i,b_i)=(0,1)$, i.e., $C_i$ is a fiber. Thus, we may have at most one fiber, at most one $-n$-curve,
	and at most one curve either in $|Z_n+nF|$ (if $n=0,1$) 
	or in $|Z_n+(n+1)F|$ (if $n=0$). Moreover, we cannot have all three.
	Thus, there are at most two components. If $r=1$ then either
	$C_1=Z_n$, so we are in case (iv), or else $C_1\in|Z_1+F|$ which is 
	case (vi), or else $C_1\in|Z_0+F|$ which is case (vii).
	If $r=2$ the only possibility is $C_1=Z_n$ and $C_2\in|F|$ (case (v)).
	
	We proceed by reductio ad absurdum to discard other cases. If $S$ had higher Picard rank, there would be a birational morphism $\pi\colon S \rightarrow s=\Fn$ for some $n\geq 0$ consisting on the consecutive contraction of $-1$-curves. Since $-K_S-C$ is ample, it follows that any $-1$-curve is either disjoint from $C$ or in the support of $C$. Since the pushforward of ample divisors is ample, we may assume without loss of generality that $\pi$ consists of one blow-up, let $c=\pi_*(C)$, $E$ be the exceptional curve, $p=\pi(E)$ and $\hat F$ be the proper transform in $C$ of the fiber $f$ through $p$ of the projection $\Fn\rightarrow \mathbb P^1$. It follows that $\hat F$ and $E$ are both $-1$-curves intersecting at one point and hence either they are both components of $C$ or both disjoint of $C$. In the latter case, it turns out that $f\cap c=\emptyset$, so $c$ consists only of fibers and by the above classification $(s,c)=(\mathbb F_0, f=Z_0)$, i.e. $(s,c)$ is of type (iv). Let $f'$ be the fiber passing through $p$ of the other projection $\mathbb F_0\rightarrow \mathbb P^1$. Then the proper transform $\hat F'$ of $F'$ via $\pi$ satisfies $(-K_S-C)\cdot \hat F'=0$, contradicting ampleness.

	Hence $\hat F,E\subseteq \mathrm{Supp}(C)$ and we will also obtain a contradiction. Write $C=\hat F+E+\Omega$ with $\hat F,E\not\subseteq\Supp(\Omega)$ and $\Omega$ some effective divisor. Let $\omega=\pi(\Omega)$. Notice that $f\not\subset \Supp(\omega)$ since $C$ is an integral divisor with coefficients at most $1$. Since $F,E\subset \Supp(C)$, we have that $-K_S-C=\pi^*(-K_S-C)+(\mult_p\omega-1)E$ and as $(-K_s-c)$ is ample, it follows that $\mult_p\omega\leq 1$. But then it follows that $(-K_S-C)\cdot E = 1-\mult_p\omega\leq 0$, contradicting ampleness.
	\epf

	\section{\Saldp surfaces}
	\lb{saldpsec}
	
	The key ideas for the classification of \saldps are quite different
	from the Maeda case. In fact, for Maeda the classification problem
	never leaves the realm of Picard rank at most $2$. In other words,
	there is no need to consider blow-ups when carrying out the classification. However, for \saldps there can be
	arbitrarily many blow-ups and there is no upper bound on the 
	Picard rank, unlike the log del Pezzo or del Pezzo settings.
	
	Thus, some sort of inductive reduction procedure is needed here.	
	The key notion is that of minimality \cite[Definition 2.8]{CR}:
	
	\bdefin
	\label{def:minimal}
	Let $(S,C)$ be \aldpno. We say it is minimal if it contains
	no smooth rational curves $E\not\subset C$ 
	with $E^2=-1$ and $E.C=1$.
	\edefin
	
	In an \aldp $(S,C)$, the dual graph of the boundary $C$ can only be either a union of chains (in fact, at most two chains \cite[Remark 3.7]{CR})
	or one cycle (in which case $C\sim -K_S$) \cite[Lemma 3.5]{CR}. In the former case we can talk of the `tail' components of each chain and the `middle' components. In the latter case all components are middle ones.
	
	Moreover, in a \saldp surface there can be no `middle' curves with negative self-intersection
	in the boundary $C$ \cite[Lemma 3.6]{CR}. Thus,
	at the very worst there are a few such curves in the `tail'. In addition, there
	are the $-1$-curves disjoint from the boundary or intersecting
	the boundary transversally at a single point.
	Contracting the latter kind, one ends up with a minimal \saldp
	which in fact has Picard rank at most 2 (this is not obvious \cite[Lemma 3.13]{CR}, see also
	\cite{R20b}) and those are readily classified. This yields
	a classification since, and this is the final step, one can verify
	generality conditions on the location of the blow-up points.
	Post factum one also finds that every pair that is obtained either as a tail
	blow-up or a blow-up away from the boundary can also be obtained by 
	blowing-up points on the boundary on a possibly different base pair. 
	Thus, tail blow-ups and `away' blow-ups play no role in the original classification of \cite{CR}. In fact, we will use some variant of this argument in the proof of Proposition \ref{partialclassifY} below.
	
	\bremark
	A systematic study of tail blow-ups has been initiated by one of us
	in the \aldp regime \cite{R19}. A systematic study of away blow-ups can be found in \cite{R20}.
	An alternative proof of the classification result of
	\cite{CR} has been given in \cite{R20b}.
	\eremark

	\section{Towards a classification
		of \aldp pairs}
	\lb{aldpsec}
	
	The key property from the classification of the strong regime
	that fails in the non-strong regime is precisely:
	no interior boundary curve with
	negative self-intersection. 
	In particular, tail blow-ups occur
	at most once in the strong regime, and one can never blow-up
	the singular points of the boundary $C$ let alone infinitely near such points.
	
	But, if one carries out the inductive proof as in the strong regime
	allowing in addition to blow down (possibly many) boundary components---which
	is never needed in the strong regime---then,
	as we will see below,  every \aldp will
	be obtained from a \aldp with Picard rank at most 2
	via {\it both} proper boundary blow-ups of smooth boundary
	points, and total boundary blow-ups of singular boundary points. However, the order of the blow-ups
	matters: one may repeatedly blow-up singular points (i.e., 
	blow-up infinitely near points) and then blow-up smooth points
	on the resulting new components of the boundary (exceptional curves).
	From the analysis of \cite{R19} it becomes clear that a daunting
	linear programming problem arises and we do not attempt to resolve it
	here.
	
	\bprop
	\lb{partialclassifY}
	Any \aldp $(S,C)$ is obtained from a pair $(s,c)$ with 
	$\rk\Pic(s)\le2$
	listed in 
	Proposition \ref{prop:classification-rk2-anticanonical} by
	a combination of the following operations:
	(i) blowing-up a collection of distinct  points on 
	the smooth locus of the boundary and replacing the boundary with its proper transform,
	(ii) blowing-up a collection of (possibly infinitely near) singular 
	points of the boundary and replacing the boundary with its total transform.
	\eprop
	
	\bremark
	Again, we emphasize that Proposition \ref{partialclassifY}
	does {\it not} give a full classification by any means,
	and it does not imply that any such blow-ups will result 
	in an \aldp pair; it merely says that any \aldp pair
	{\it can be described in such a way}, but we do not give
	that description here.
	\eremark

	\bpf
	Reversing both operations (i) and (ii) preserves asymptotic log positivity.
	Indeed, contraction of a $-1$-curve intersecting the boundary transversally at a single point preserves the property of being \aldp \cite[Lemma 3.4]{CR}.
	And if
	$C=C_1+E+\sum_{i=3}^rC_i$
	with $E=C_2$ being a $-1$-curve,
	and $E.C_i=1$ for $i=1,3$ and zero otherwise with blow-down $\pi:S\ra s$ contracting $E$
	to a point, then letting $c_1:=\pi(C_{1}), c_i:=\pi(C_{i+1}), i=2,\ldots,r-1$,
	$$
	-K_S-\sum_{i=1}^r(1-\be_i)C_i
	\sim
	-\pi^*\Big(K_s+\sum_{i=1}^{r-1}(1-\be_i)c_i\Big)-(\be_1+\be_3-\be_2)E,
	$$
	so 
	$$
	-K_s-\sum_{i=1}^{r-1}(1-\be_i)c_i\sim\pi_*
	\Big(-K_S-\sum_{i=1}^r(1-\be_i)C_i\Big)
	$$ is therefore ample by the Nakai--Moishezon criterion (see, e.g., \cite[(1.5)]{CR}),
	so $(s,\sum_{i=1}^{r-1}c_i)$ is \aldpno. 
	
	Next, we claim that we can reverse the operations (i) and (ii)
	(in some order), until
	the Picard rank becomes at most 2 (while preserving
	the \aldp property as we just showed above).
	To prove the claim, first apply operations (i)
	to contract all $-1$-curves in $S$ that intersect the boundary 
	transversally at a smooth point (note that these are in fact all $-1$-curves
	not contained in $C$ that intersect the boundary \cite[Lemma 3.3]{CR}).
	Second, apply operations (ii) to contract any remaining
	$-1$-curve in the resulting boundary
	that intersects two distinct other boundary components (when $r\ge3$).
	We can also assume that $r\ge2$ since if $r=1$ the pair is \saldp 
	and so by the classification \cite[Theorem 2.1]{CR} we can contract it down
	using only the operations (i). 
	If after the above operations the the rank of the Picard group of the 
	resulting surface is at most 2 we are done. 
	Thus, we may assume that we are in one of the following mutually exclusive situations:
	(a) $\rk\Pic(S)\ge 3,\, r\ge2,\, C\not\sim-K_S$, $C$ contains some $-1$-tail components
	and the only other $-1$-curves in $S$ are disjoint from $C$, 
	or 
	(b) $\rk\Pic(S)\ge 3,\, r\ge2,$, the only $-1$-curves in $S$ are disjoint from $C$,
	or
	(c) $\rk\Pic(S)\ge 3,\, r=2,\, C_1+C_2\sim-K_S$ with $C_1^2=-1$ and $C_1.C_2=2$.
	Here we used \cite[Lemma 3.5]{CR} for the structure of $C$.
	
	To treat (a), contract first all of the $-1$-curves disjoint from $C$.  The resulting pair,
	that we will still denote by $(S,C)$ (without loss of generality) is \aldp \cite[Lemma 3.4]{CR}. Furthermore, 
	this contraction still preserves the fact that there are no $-1$-curves not contained in the boundary that intersect the remaining tails: the contraction can only increase the self-intersection of 
	curves in $S$ but it is impossible to have a curve of negative self-intersection less than $-1$ not in $C$ to begin with
	\cite[Lemmas 2.5, 2.7]{CR}. This contraction also evidently preserves the property that 
	there are still no $-1$-curves in $C$ that intersect two other components of $C$. 
	Next, contract one of the $-1$-tails in $C$. The new boundary will have the property that it contains
	no `middle' $-1$-curves while if the boundary component that intersected that tail upstairs was a $-2$-curve
	it will now be itself a $-1$-tail. We may therefore repeat this process 
	(this corresponds to reverse operation to blowing-up infinitely near points)
	until the rank of the Picard group is exactly 3. By induction, in this process there is never a `middle' 
	$-1$-component of the boundary.
	Denote by $(S,C)$ (without loss of generality) the resulting 
	\aldp pair \cite[Lemma 3.12]{CR} with $\rk\Pic(S)=3$. 
	After this most recent contraction there is still at least one remaining $-1$-tail, say $C_1\subset C$, 
	but no `middle' $-1$-curves in $C$
	(since every rational surface with Picard group of rank 3 has a $-1$-curve and we know the only such curves
	in $S$ must be tail components of $C$ by our careful construction),	and furthermore
	that $C_1$ still does not intersect any $-1$-curves not contained in the boundary.
	Recall that rational surfaces of Picard rank $3$ have precisely three curves of negative self-intersection namely $-k$, $-1$, $-1$, $k>0$, with a $-1$-curve being the `middle' curve in such chain. The curve $C$ must contain one of these $-1$-curves, and by assumption it cannot contain a `middle' $-1$-curve and it cannot intersect a $-1$-curve not in $C$, so all $-1$-curves are in $C$. If $k=1$ this means $C$ has a `middle' $-1$-curve, giving a contradiction. If $k>1$, the $-k$-curve must be in $C$ by \cite[Lemma 2.5]{CR}, giving a $-1$-curve as a `middle' curve. In conclusion then we obtain a contradiction which concludes our treatment of (a).
	
	To treat (b), note that $C\not\sim-K_S$ as for any $-1$-curve $E$ have $-K_S.E=1$ by adjunction 
	which forces $E.C_i=1$ for some $i$ if $E\not\subset C$. Contract some of the $-1$-curves on $S$ until
	the Picard group has rank 3. As in case (a) this process still preserves the property that 
	there are no $-1$-curves not in the boundary that intersect the boundary. As in (a), there are $3$-curves of negative self-intersection forming a chain. By \cite[Lemma 2.5]{CR}, any curve $Z$ with $Z^2<-1$ must be in $C$ and thus intersect a $-1$-currve, forcing $S$ to be the del Pezzo surface of degree $7$, with a chain of three $-1$-curves $e_1, l, e_2$, not intersecting $C$. We have that $-K_S\sim 2e_1+2e_2+3l>0$ so one of these curves must intersect $C$ by ampleness, giving a contradiction and concluding our treatment of (b).

	To treat (c), observe that $S$ contains no $-1$-curves $E$ disjoint from $C$ since
	$-K_S.E=1$ by adjunction, but on the other hand $-K_S\sim C_1+C_2$ implies either $E.C_1=1$ or $E.C_2=1$. By our
	construction $S$ also contains no $-1$-curves intersecting $C$ transversally so it follows that all $-1$-curves
	on $S$ are in $C$. If $C_1$ is the only $-1$-curve in $S$ then $\rk\Pic(S)=2$ and we are done. So assume that
	$\rk\Pic(S)=3$ and that $C_1$ and $C_2$ are all the $-1$-curves on $S$. Again, by the classification
	of rational surfaces the only rational surfaces with exactly two $-1$-curves and Picard group of rank 3 are the blow-up of $\FF_n,\,
	n\ge2$ at a point or the blow-up of $\FF_1$ on a point in $Z_1$. 
	However, in all of these cases the two $-1$-curves intersect
	transversally at a single point and not in two points as we assumed in (c). Thus, (c) cannot happen and the proof of Proposition \ref{partialclassifY} is complete.
	\epf

	\bremark
	Proposition \ref{partialclassifY} shows that \aldp pairs are a subset of a rather explicit
	(infinite) family of pairs. Thus, it reduces Problem \ref{mainprob} to determining the generality conditions
	on the blow-ups of type (i) and (ii) of the pairs listed in Proposition \ref{prop:classification-rk2-anticanonical}.
	As an example, 
	if $(s,c)$ is 
		$(\mathrm{ALdP.4.n})$ and if $f\not\subset c$ is a fiber, we may not blow-up two smooth points on $f\cap c$
	since then for any choice of $\beta_i$ arbitrarily close to $0$, 
	$(-K_S-\sum_{i=1}^4(1-\be_i)C_i).\tilde f=0$.
	Hence, 
	one generality condition would be ``no two of the points blown-up may lie on a fiber not contained in the boundary." 
	Part of the difficulty in describing such generality conditions lies on the fact that there can infinitely-near points of any order in the operations (ii) and so one needs
	to describe these generality conditions on all the `intermediate surfaces', so to speak.
	\eremark

	\section{Classification of \aldp
		pairs with small Picard group}
	\lb{aldppic2sec}
	
	In this section we classify \aldps $(S,C)$ with  $\rk(\Pic(S))\le 2$.

	\bprop
	\label{prop:classification-rk2-anticanonical}
	Let $S$ be a smooth surface with $\rk(\Pic(S))\le 2$, and let
	$C_1,\ldots,C_r$ be distinct irreducible smooth curves on $S$ such that
	$C=\sum_{i=1}^{r}C_i$  is a divisor with simple normal crossings.
	Then $(S,C)$ is an 
	asymptotically log del Pezzo pair if and only if it
	is one of the following pairs:
	\begin{itemize}
		\item [$\mathrm{(I.1A})$] $S=\mathbb{P}^2$, $C_1$ is a cubic,%
		\item [$\mathrm{(I.1B})$] $S=\mathbb{P}^2$, $C_1$ is a conic,%
		\item [$\mathrm{(I.1C})$] $S=\mathbb{P}^2$, $C_1$ is a line,%
		\item [$\mathrm{(I.2.n})$] $S=\mathbb{F}_{n}, \,\,n\in\NN\cup\{0\}$, $C_1=Z_n$,%
		\item [$\mathrm{(I.3A})$] $S=\mathbb{F}_1$, $C_1\in |2(Z_1+F)|$,%
		\item [$\mathrm{(I.3B})$] $S=\mathbb{F}_1$, $C_1\in |Z_1+F|$,%
		\item [$\mathrm{(I.4A})$] $S=\mathbb{P}^1\times\mathbb{P}^1$,  $C_1$ is a $(2,2)$-curve,%
		\item [$\mathrm{(I.4B})$] $S=\mathbb{P}^1\times\mathbb{P}^1$,  $C_1$ is a $(2,1)$-curve,%
		\item [$\mathrm{(I.4C})$] $S=\mathbb{P}^1\times\mathbb{P}^1$,  $C_1$ is a $(1,1)$-curve,%
		\item [$\mathrm{(I.5.1})$] $S=\mathbb F_1$, $C_1\in|2Z_1+3F|$, 
		\item [$\mathrm{(I.6B.1})$] $S=\mathbb F_1$, $C_1\in |Z_1+2F|$,
		\item [$\mathrm{(I.6C.1})$] $S=\mathbb F_1, C_1\in |F|$,
		\item [$\mathrm{(II.1A})$]  $S=\mathbb{P}^2$, $C_1$ is a conic, $C_2$ is a line, 
				\item [$\mathrm{(II.1B})$]  $S=\mathbb{P}^2$, $C_1,C_2$ are lines,%
				\item [$\mathrm{(II.2A.n})$]  $S=\mathbb{F}_n, \,n\in\NN\cup\{0\}$,  $C_1=Z_n, C_2\in |Z_n+nF|$,%
				\item [$\mathrm{(II.2B.n})$]  $S=\mathbb{F}_n, \,n\in\NN\cup\{0\}$,  $C_1=Z_n, C_2\in|Z_n+(n+1)F|$,%
				\item [$\mathrm{(II.2C.n})$]  $S=\mathbb{F}_n, \,n\in\NN\cup\{0\}$, $C_1=Z_n, C_2\in|F|$,%
				\item [$\mathrm{(II.3})$]  $S=\mathbb{F}_1$, $C_1, C_2\in|Z_1+F|$,
				\item [$\mathrm{(II.4A})$]  $S=\mathbb{P}^1\times\mathbb{P}^1$,  $C_1,C_2$ are $(1,1)$-curves,
				\item [$\mathrm{(II.4B})$]  $S=\mathbb{P}^1\times\mathbb{P}^1$, 
		$C_1$ is a $(2,1)$-curve, $C_2$ is a $(0,1)$-curve, 
				\item [$\mathrm{(II.5A.1})$] $S=\mathbb F_1$, either $C_1\in |2Z_1+2F|, C_2\in |F|$,
		and $|C_1\cap C_2|=2$, or else $C_1\in |Z_1+2F|, C_2\in |Z_1+F|$, 
		\item [$\mathrm{(II.5B.1})$] $S=\mathbb F_1$, $C_1\in |Z_1+F|, C_2\in |F|$,
		\item [$(\mathrm{ALdP.1.n})$
		]  $S= \mathbb F_n, \, n\in\NN$, $C_1=Z_n, C_2\in |Z_n+(n+2)F|$,
		
		\item [$\mathrm{(III.1})$]  $S=\mathbb{P}^2$, $C_1, C_2, C_3$ are lines, 
		
		\item [$\mathrm{(III.2})$]  $S=\mathbb{P}^1\times\mathbb P^1$, 
		$C_1, C_2, C_3$ are $(1,1)$-, $(0,1)$-, and $(1,0)$-curves, respectively,
		
		\item [$\mathrm{(III.3.n})$]  $S=\mathbb{F}_{n}, \,\,n\in\NN\cup\{0\},  \, C_1=Z_n, C_2\in |F|, 
		C_3\in |Z_n+nF|$,
		
		\item [$\mathrm{(III.4.1})$]  $S=\mathbb{F}_1$, $C_1\in |F|, C_2, C_3\in |Z_1+F|$,
		
		\item [$(\mathrm{ALdP.2.n})$
		] $S= \mathbb F_n, \, n\in\NN, \, C_1=Z_n, C_2\in |Z_n+(n+1)F|, C_3\in|F|$,
		
		\item [$(\mathrm{ALdP.3.n})$
		] $S= \mathbb F_n, \, n\in\NN, \,C_1=Z_n, C_2, C_3\in |F|$,
		
		\item [$\mathrm{(IV})$]
		$S=\mathbb{P}^1\times\mathbb{P}^1$, $C_1, C_2$ 	are $(1,0)$-curves, $C_3,C_4$ are $(0,1)$-curves,
		
		\item [$(\mathrm{ALdP.4.n})$
		] $S= \mathbb F_n, \, n\in\NN, \,C_1=Z_n, C_2, C_3\in|F|,, C_4\in |Z_n+nF|$.
		
	\end{itemize}
	
	\eprop
	
	\bremark
	Note that of course when we say, e.g., ``$C_1$ is a conic",
	we mean it is a smooth conic since we assume, as in the statement
	that each $C_i$ is smooth. Similarly, e.g., in (II.1A) the line and the
	conic {\it must} meet at two distinct point as we require $C_1+C_2$
	to have simple normal crossings.
	\eremark
	
	\begin{proof}
		We start with the sufficiency statement. 
		All pairs in the statement, except (ALdP.1.n),$\ldots,$ (ALdP.4.n), are \saldp	\cite[Theorems 2.1,3.1]{CR}. The remaining four cases are \aldp but not \saldp as can be seen from 
		\er{ampleFn} which also computes their bodies
		of ample angles in terms of a single inequality, that we have recorded in Corollary 
		\ref{aldpnotsaldpcor} below and in Figures \ref{Figure1}--\ref{AA-2}.
		
		Let us turn to the necessity statement, i.e., suppose that $(S,C)$ is \aldp with
		$\rk\Pic(S)\le 2$.
		First, $S$ is rational, hence it is either $\PP^2$ or $\Fn$
		\cite[\S2--3]{CR}. When $S=\PP^2$, $\rk(\Pic(S))=1$ so every \aldp is automatically \saldpno,
		and since $-K_S\sim 3H$ we see the possibilities are 
		$\mathrm{(I.1A}),\mathrm{(I.1B}),\mathrm{(I.1C})$, $\mathrm{(II.1A}),\mathrm{(II.1B}),\mathrm{(III.1})$.
		Assume now that $S=\Fn$. Denote
		$$
		C_i\in|a_iZ_n+b_iF|.
		$$
		Since $-K_S-C$ is nef by \cite[\S2.1]{CR}  and using \er{KFnEq} and \er{nefFn}, we see that 
		$$
		\sum_i a_i\in\{0,1,2\}, \q 
		\sum_i b_i\in\{0,\ldots,n+2\}. 
		$$
		Note that by \er{irreducibleFn} either $b_i\ge na_i>0$, or $(a_i, b_i)=(1,0)$, i.e. $C_i=Z_n$ in which case this can only happen for one $i$, or else $(a_i,b_i)=(0,1)$, i.e., $C_i$ is a fiber.
		Thus, at most two components of $C_i$ are {\it not } fibers. Since every fiber intersects any other curve that is not a fiber by at least 1 and by \cite[Lemma 3.5]{CR} the dual graph of $C$ is either a union of chains or a cycle, if one $C_i$ is not a fiber, then there are at most two fibers in $C$. In particular $1\le r\le 4$ or $C$ consists only of fibers. Also, if $a_i=2$ for some $i$ 
		then $n\le 2$ by $n+2\ge b_i\ge 2n$. So we get the following possibilities when
		$\max_ia_i=2$:
		\beq
		\lb{eqcases1}
		\baeq
		{}[n,(a_1,b_1),\ldots,(a_r,b_r)]\in\big\{
		&[2,(2,4)],[1,(2,3)],[1,(2,2),(0,1)],
		\cr
		&[1,(2,2)],[0,(2,2)],[0,(2,1),(0,1)],
		[0,(2,1)]\big\}.
		\eaeq
		\eeq
		When $\max_ia_i=1$, we split into two subcases:
		when
		there are at least two pairs $(a_i,b_i)$ with all coefficients
		positive (in which case $n+2\ge \sum_ib_i\ge 2n$, as $b_i\ge na_i=n$ for at least two coefficients, so again $0\le n \le 2$):
		\beq
		\lb{eqcases2}
		\baeq
		{}[n,(a_1,b_1),\ldots,(a_r,b_r)]\in\big\{
		&[2,(1,2),(1,2)],[1,(1,2),(1,1)],
		\cr
		&[1,(1,1),(1,1),(0,1)],
		[1,(1,1),(1,1)],[0,(1,1),(1,1)]
		\big\},
		\eaeq
		\eeq
		and otherwise, still with $\max_ia_i=1$, now for all $n$ (so we omit the first index), and
		with $\max_ib_i=n$,
		\beq
		\lb{eqcases3}
		\baeq
		{}[(a_1,b_1),\ldots,(a_r,b_r)]\in\big\{
		&[(1,n),(1,0),(0,1),(0,1)],
		[(1,n),(1,0),(0,1)],
		\cr
		&[(1,n),(1,0)],
		[(1,n),(0,1),(0,1)],
		[(1,n),(0,1)],
		[(1,n)]
		\big\},
		\eaeq
		\eeq
		with $\max_ib_i=n+1$,
		\beq
		\lb{eqcases4}
		\baeq
		{}[(a_1,b_1),\ldots,(a_r,b_r)]\in\big\{
		&[(1,n+1),(1,0),(0,1)],
		\cr
		&[(1,n+1),(1,0)],
		[(1,n+1),(0,1)],
		[(1,n+1)]
		\big\},
		\eaeq
		\eeq
		with $\max_ib_i=n+2$,
		\beq
		\lb{eqcases5}
		\baeq
		{}[(a_1,b_1),\ldots,(a_r,b_r)]\in\big\{
		&[(1,n+2),(1,0)],
		[(1,n+2)]
		\big\},
		\eaeq
		\eeq
		when $\max_ib_i=1$,
		\beq
		\lb{eqcases6}
		\baeq
		{}[(a_1,b_1),\ldots,(a_r,b_r)]
		\in\big\{
		&[(1,0)],
		[(1,0),(0,1)],
		[(1,0),(0,1),(0,1)],
		\big\},
		\eaeq
		\eeq
		since as we saw earlier if $C$ has an irreducible curve in its support that it is not a fiber, then $C$ can have at most $2$ fibers in its support.
		
		Finally when $\max_ia_i=0$,
		\beq
		\lb{eqcases7}
		\baeq
		{}[(a_1,b_1),\ldots,(a_r,b_r)]
		\in\big\{
		&[(0,1)],\ldots,
		[\overbrace{(0,1),\ldots,(0,1)}^{n+2}]
		\big\}.
		\eaeq
		\eeq
		A few of these cases can be eliminated, though most of them
		actually occur. 
		In 
		\er{eqcases1}, [2,(2,4)]
		corresponds to a smooth anticanonical curve in $\FF_2$, 
		which is excluded as $\FF_2$ is not del Pezzo.
		The remaining cases are: $[1,(2,3)]=$  I.5.1,
		$[1,(2,2),(0,1)]=$ II.5A.1 (blow-up on the line in II.1A),
		$[1,(2,2)]=$ I.3A,
		$[0,(2,1),(0,1)]=$ II.4B,
		$[0,(2,1)]=$ I.4B.
		
		In \er{eqcases2}
		[2,(1,2),(1,2)] is excluded as $Z_2.(Z_2+2F)=Z_2.(2Z_2+4F)=Z_2.(-K_{\FF_2})=0$.
		The remaining cases are: 
		$[1,(1,2),(1,1)]=$ II.5A.1  (blow-up on the conic in II.1A);\hfill\break
		$[1,(1,1),(1,1),(0,1)]=$ III.4.1,
		$[1,(1,1),(1,1))=$ II.3,
		$(0,(1,1),(1,1)]=$ II.4A. 
		
		In \er{eqcases3},
		$[(1,n),(1,0),(0,1),(0,1)]$
		gives $-K_{\FF_n}-(1-\be_1)(Z_n+nF)
		-(1-\be_2)Z_n
		-(1-\be_3)F-(1-\be_4)F\sim(\be_1+\be_2)Z_n+(n\be_1+\be_3+\be_4)F$,
		that is ample if and only if 
		$n\be_1+\be_3+\be_4> n\be_1+n\be_2$, i.e., 
		$\be_3+\be_4>n\be_2$, and this is ALdP.4.n if $n\ge 1$ or IV if $n=0$; 
		$[(1,n),(1,0),(0,1)]=$ III.3.n,
		$[(1,n),(1,0)]=$ II.2A.n.
		For [(1,n),(0,1),(0,1)] consider 
		$-K_{\FF_n}-(1-\be_1)(Z_n+nF)
		-(1-\be_2)F-(1-\be_3)F\sim(1+\be_1)Z_n+(n\be_1+\be_2+\be_3)F$,
		that is ample if and only if 
		$n\be_1+\be_2+\be_3> n+n\be_1$, i.e., $n=0$, and this is
		III.3.0.
		For [(1,n),(0,1)], consider 
		$-K_{\FF_n}-(1-\be_1)(Z_n+nF)
		-(1-\be_2)F\sim(1+\be_1)Z_n+(1+n\be_1+\be_2)F$,
		that is ample if and only if 
		$1+n\be_1+\be_2> n+n\be_1$, i.e., $n=0,1$, and these are
		II.2C.0, II.5B.1.
		For [(1,n)], $-K_{\FF_n}-(1-\be_1)(Z_n+nF)\sim
		(1+\be_1)Z_n+(2+n\be_1)F$ implying $n=0,1,2$ and
		these are I.2.0, I.3B, while the case $n=2$ is excluded
		as in the first paragraph.
		
		In \er{eqcases4}:
		[(1,n+1),(1,0),(0,1)] 
		$-K_{\FF_n}-(1-\be_1)(Z_n+(n+1)F)
		-(1-\be_2)Z_n-(1-\be_3)F\sim(\be_1+\be_2)Z_n+((n+1)\be_1+\be_3)F$
		that is ample if and only if
		$(n+1)\be_1+\be_3>n(\be_1+\be_2)$, i.e., $\be_1+\be_3>n\be_2$,
		and this is ALdP.2.n if $n\ge 1$ and III.2 if $n=0$;
		$[(1,n+1),(1,0)]=$ II.2B.n;
		[(1,n+1),(0,1)],
		$-K_{\FF_n}-(1-\be_1)(Z_n+(n+1)F)
		-(1-\be_2)F\sim(1+\be_1)Z_n+((n+1)\be_1+\be_2)F$
		that is ample if and only if
		$(n+1)\be_1+\be_2>n(1+\be_1)$, i.e., $n=0$
		and this is II.2B.0;
		[(1,n+1)],
		$-K_{\FF_n}-(1-\be_1)(Z_n+(n+1)F)
		\sim(1+\be_1)Z_n+(1+(n+1)\be_1)F$
		that is ample if and only if
		$1+(n+1)\be_1>n(1+\be_1)$, i.e., $n=0,1$,
		and these are I.4C, I.6B.1.
		
		In \er{eqcases5},
		[(1,n+2),(1,0)],
		$-K_{\FF_n}-(1-\be_1)(Z_n+(n+2)F)
		-(1-\be_2)Z_n\sim(\be_1+\be_2)Z_n+(n+2)\be_1F$,
		i.e., $2\be_1>n\be_2$, and this is ALdP.1.n; 
		[(1,n+2)],
		$-K_{\FF_n}-(1-\be_1)(Z_n+(n+2)F)
		\sim(1+\be_1)Z_n+(n+2)\be_1F$,
		i.e., $(n+2)\be_1>n+n\be_1$, i.e., $n=0$ and
		this is I.4B.
		
		In \er{eqcases6},
		there is one $Z_n$ and $0\le k\le 2$ fibers. When $k=0$ this is I.2.n, when $k=1$ this is
		II.2C.n, and when $k=2$:
		$-K_{\FF_n}-(1-\be_1)Z_n
		-(1-\be_2)F-
		(1-\be_3)F
		\sim(1+\be_1)Z_n+(n+\be_2+\be_3)F$,
		that is ample if and only if 
		$
		n(1+\be_1)
		<
		(n+\be_2+\be3)
		$, and is
		is ALdP.3.n. 
		
		Finally, in \er{eqcases7},
		there are $k$ fibers, so
		$-K_{\FF_n}-(1-\be_1)F
		-\ldots
		(1-\be_{k})F
		\sim2Z_n+(n+2-k+\be_1+\ldots+\be_{k})F$,
		that is ample if and only if 
		$
		2n
		<
		(n+2-k+\be_1+\ldots+\be_{k})
		$,
		i.e., $n=k=1$ and I.6C.1, or $n=0$ and $k=1,2$
		and these are I.2.0, II.2A.0.
		\epf
		
		\bcor
		\lb{aldpnotsaldpcor}
		Let $S$ be a smooth surface with $\rk(\Pic(S))\le 2$, and let
		$C_1,\ldots,C_r$ be irreducible smooth curves on $S$ such that
		$C=\sum_{i=1}^{r}C_i$  is a divisor with simple normal crossings.
		Then $(S,C)$ is an 
		asymptotically log del Pezzo pair, but not \saldpno, if and only if it
		is one of the pairs \emph{ALdP.1.n, ALdP.2.n, ALdP.3.n, ALdP.4.n.}
		Moreover,
		\begin{equation*}
		\h{\rm AA}(S,C)=
		\begin{cases}
		\{(\be_1,\be_2)\in(0,1]^2\,:\,-n\be_1+2\be_2>0\}
		\h{ if  $(S,C)$ is $\mathrm{ALdP.1.n}$}, 
		\cr
		\{(\be_1,\be_2,\be_3)\in(0,1]^3\,:\,-n\be_1+\be_2+\be_3>0\}
		\h{ if  $(S,C)$ is $\mathrm{ALdP.2.n}$ or $\mathrm{ALdP.3.n}$}, \cr
		\{(\be_1,\be_2,\be_3,\be_4)\in(0,1]^4\,:\,-n\be_1+\be_2+\be_3>0\}
		\h{ if  $(S,C)$ is $\mathrm{ALdP.4.n}$}. 
		\end{cases}
		\end{equation*}
		\ecor

		\begin{figure}
			\centering
			\includegraphics[width=0.3\textwidth]{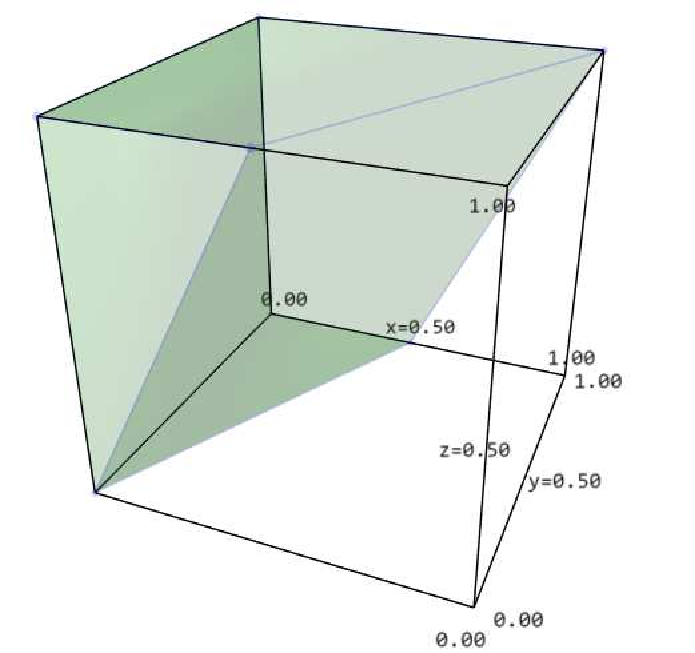}
			\caption{The body of ample angles for $\mathrm{ALdP.k.2}$ for $k=2,3$ and the projection of $\mathrm{ALdP.4.2}$ on $\beta_1=m$ hyperplane for $m\in (0,1)$.}
			\label{AA-2}
		\end{figure}
		
		\section{The body of ample angles}
		\lb{polytopesec}
		
		The goal of this section is to prove Theorem \ref{CasciniRThm}, as a consequence of the main result in the theory of Shokurov's log geography (cf.  \cite{BCHM,CL,KKL,SC}). 
		
		Before we proceed with the proof, we recall some of the  notions and results which we will use later. As in Section 1,  $(X,D)$ is a pair consisting of a smooth complex projective variety $X$ and a simple-normal crossing divisor 
		$D=\sum_{i=1}^rD_i$. Let $\rm{WDiv}_{\mathbb R}(X)$ be the real vector space spanned by all the divisors in $X$, and let $V\subset \rm{WDiv}_{\mathbb R}(X)$ be the finite-dimensional vector subspace spanned by the basis $D_1,\dots,D_r$, and use that basis to identify $V\cong \mathbb R^r$. Fix an ample $\mathbb Q$-divisor $A$ and define 
		$$\mathcal L(V):=\{B\in V: B=\sum_{i=1}^r a_i D_i\quad\text{where }a_1,\dots,a_r\in [0,1]^r \},$$ 
		and 
		$$\mathcal A_A(V):=\{B\in \mathcal L(V): K_X+A+B \text{ is ample}\}.$$
		
		We want to show that the set $\mathcal A_A(V)$  is a rational polytope. Indeed, let $\psi\colon X\to X$ be the identity morphism and let $\pi\colon X\to \rm{Spec}\, \mathbb C$ be the structure morphism.
		By definition, $\rm{Spec}\, \mathbb C$ is a point, so 
		$\pi$ is the trivial morphism to a point.
		Then, $\mathcal A_A(V)$ coincides with $A_{\psi,A,\pi}(V)$, as defined in \cite[Definition 1.1.4]{BCHM}.
		Thus, by \cite[Corollary 1.1.5]{BCHM}, it follows that $\mathcal A_A(V)$ is a finite union of rational polytopes and, by convexity of $\mathcal A_A(V)$, the claim follows. 
		
		We now proceed with the proof of our Theorem. 
		
		\begin{proof}[Proof of Theorem \ref{CasciniRThm}]
			Suppose $\hbox{\rm AA}(X,D)$ is non-empty. 
			We claim that the closure of $\hbox{\rm AA}(X,D)$ is a
			rational polyhedron.
			By openness of ampleness, there exists  $\gamma=(\gamma_1,\dots,\gamma_r)\in (0,1]^r\cap\QQ$ such that $-K_X-\sum_{i=1}^r(1-\gamma_i)D_i$ is ample. 
			We fix such a $\gamma$ for the rest of the proof.
			Let 
			$$
			\eta:=\max\bigg\{\max_i\frac{1-\gamma_i}{\gamma_i},
			\max_i\frac{\gamma_i}{1-\gamma_i}\bigg\}.
			$$
			For any $\beta\in(0,1)^r$ we claim we can express
			\beq
			\lb{etaFormula}
			-K_X-\sum_{i=1}^r(1-\beta_i)D_i
			=
			\eta
			\big(
			K_X+A+F(\be)
			\big),
			\eeq
			with two key properties:
			
			\smallskip	
			\noi $\bullet\;$
			$A$ ample and independent of $\be$, and 
			
			\smallskip	
			\noi 
			$\bullet\;$
			$F(\be)$ effective, supported on $D$, and with the coefficient
			of each $D_i$ in $(0,1)$ (and, of course,
			depending on $\be$), in such a way that the map between
			$\be\in(0,1)^r$ and the vector of coefficients of $F(\be)$ 
			in $(0,1)^r$ is one-to-one.
			
			\smallskip	
			More specifically, let
			$$
			A:=-\frac{1+\eta}{\eta}\bigg(K_X+\sum(1-\gamma_i)D_i\bigg),
			$$
			and
			$$
			F(\be):=
			\sum_{i=1}^r\bigg(1+\frac{\beta_i}\eta - \frac{1+\eta}{\eta}\gamma_i\bigg)D_i
			.$$
			Then, one readily checks that \er{etaFormula} holds for
			each $\be\in(0,1)^r$. It remains to check the two key properties.
			First, $A$ is ample by definition since $\eta>0$ and $\gamma\in \hbox{\rm AA}(X,D)$;
			$A$ is also clearly independent of $\be$. 
			Second, we need to verify that
			\beq
			\lb{etaIneq}
			0\le 
			1+\frac{\beta_i}\eta - \frac{1+\eta}{\eta}\gamma_i
			\le 1,
			\eeq
			for each $i=1,\ldots,r$ and each $\be\in(0,1)^r$.
			The upper bound amounts to
			$$
			\frac{\beta_i}\eta \le \frac{1+\eta}{\eta}\gamma_i
			,
			$$
			i.e.,
			$$
			\frac{\beta_i}{\gamma_i}-1 \le\eta
			,
			$$
			and this holds since 
			$$
			\frac{\beta_i}{\gamma_i}-1 \le
			\frac{1-\gamma_i}{\gamma_i}\le
			\max_i\frac{1-\gamma_i}{\gamma_i}\le
			\eta.
			$$
			The lower bounds amounts to
			$$
			0\le 
			1+\frac{\beta_i}\eta - \frac{1+\eta}{\eta}\gamma_i,
			$$
			i.e.,
			$$
			\frac{\gamma_i-\beta_i}{1-\gamma_i}\le\eta
			,
			$$
			and this holds since 
			$$
			\frac{\gamma_i-\beta_i}{1-\gamma_i}
			\le
			\frac{\gamma_i}{1-\gamma_i}\le
			\max_i\frac{\gamma_i}{1-\gamma_i}\le
			\eta.
			$$
			In sum, we have proved the following. Let
			$
			f:\RR^r\ra\RR^r
			$ 
			be the affine map
			$$
			f(\be):=
			(1,\ldots,1)+\frac{\beta}\eta - \frac{1+\eta}{\eta}\gamma.
			$$
			Then $f$ maps the cube $(0,1)^r$ into $(0,1)^r$.
			We claim that these facts combined imply that
			\beq
			\lb{beautiful}
			f^{-1}(\mathcal A_A(V))\cap (0,1)^r=\hbox{\rm AA}(X,D),
			\eeq
			where we identify $\mathcal A_A(V))$ with a subset of $(0,1)^r$ via the basis $D_1, \ldots, D_r$.
			Indeed, if $\be\in \hbox{\rm AA}(X,D)$ then 
			by definition $\be\in(0,1)^r$, and, moreover, by
			\er{etaFormula} $f(\be)\in \mathcal A_A(V)$, i.e.,
			$\be\in f^{-1}(\mathcal A_A(V))\cap (0,1)^r$. Conversely,
			if 
			$\be\in f^{-1}(\mathcal A_A(V))\cap (0,1)^r$, write
			$\be=f^{-1}(\al)$ with $\al \in \mathcal A_A(V)$.
			One readily checks that $f^{-1}$ can be written
			explicitly as
			\beq\lb{finv}
			f^{-1}(\al):=
			\eta\al-(\eta,\ldots,\eta)
			+(1+\eta)\gamma.
			\eeq
			By definition of $\al$ we know that
			$
			K_X+A+\sum\al_iD_i,
			$
			and hence also
			$
			\eta(K_X+A+\sum\al_iD_i),
			$
			is ample. But the latter equals identically
			$-K_X-D+\eta\sum\al_iD_i+(1+\eta)\sum\gamma_iD_i-\eta D$
			which by \er{finv} equals 
			$-K_X-D+\sum (f^{-1}(\al))_iD_i=-K_X-\sum(1-\be_i)D_i$;
			recalling that by assumption $\be\in(0,1)^r$,
			we have shown that 
			$\be\in\hbox{\rm AA}(X,D)$. Altogether,
			we have established \er{beautiful}, as desired.
			We remark in passing that it
			is not true in general that
			$f^{-1}$ maps the cube $(0,1)^r$ into itself, but
			that is not needed for our proof.
			Finally, using \er{beautiful}, that
			$\mathcal A_A(V)$ is a rational polyhedron by the discussion at the beginning of this section, and since $f$ is an affine
			map, it follows that $\hbox{\rm AA}(X,D)$ is a rational polyhedron.
		\end{proof}
		
		\bremark We give a second proof of Theorem \ref{CasciniRThm} that is a bit more highbrow.
		Let $X$ be a smooth projective variety. We denote by $NS(X)$ the N\'eron-Severi group of $X$, i.e. the  real vector space spanned by the divisors of $X$ modulo numerical equivalence.  We denote by ${\rm Nef}(X)\subset NS(X)$ the { nef cone} of $X$, i.e., the cone spanned by the
		numerically effective (nef) divisors in $NS(X)$. 
		Let $L_1,\dots,L_k$ be divisors in $X$. The  {\it Cox ring} of $X$ associated to $L_1,\dots,L_k$ is the ring 
		\[R(X,L_\bullet)\coloneqq \bigoplus_{m_1,\dots,m_k\ge 0}H^0\left(X,\mathcal O_X\left(\sum_{i=1}^k mL_i\right)\right). 
		\]
		The variety $X$ is called a {\it Mori dream space} if $H^1(X,\mathcal O_X)=0$ and there exist divisors $L_1,\dots,L_k$ whose numerical classes generate $NS(X)$ and such that the ring $R(X,L_\bullet)$ is finitely generated. 
		If $X$ is a Mori dream space then, in particular, ${\rm Nef}(X)$ is a rational polyhedral (cf. \cite[Proposition 2.9 and Definition 1.10]{HK00}). 
		Assume now that $D=\sum_{i=1}^r D_i$ is a simple-normal-crossing divisor in $X$ such that  $\hbox{\rm AA}(X,D)$ is not empty and let $\gamma=(\gamma_1,\dots,\gamma_r)\in\hbox{\rm AA}(X,D)$. Let \[
		\Delta=\sum_{i=i}^r(1-\gamma_i)D_i. 
		\]
		Then $(X,\Delta)$ is log smooth and, in particular, it is divisorially log terminal (e.g. see \cite[Definition 3.1.1]{BCHM}). Thus, \cite[Corollary 1.3.2]{BCHM} implies that $X$ is a Mori dream space.
		Consider the affine map 
		\[
		\Phi\colon\mathbb R^r\to NS(X)\]
		given by, for any $(\beta_1,\dots,\beta_r)\in \mathbb R^r$, 
		\[\Phi(\beta_1,\dots,\beta_r)= [-K_X-\sum_{i=1}^r(1-\beta_i)D_i].
		\]
		Then, from the definitions, 
		\[
		\overline{\hbox{\rm AA}(X,D)}=[0,1]^r\cap \Phi^{-1}({\rm Nef}(X))
		\]
		(here we denote by $\Phi^{-1}({\rm Nef}(X))$ the pre-image of 
		${\rm Nef}(X)$).
		Since $NS(X)$ is a finite dimensional vector space, we may fix an isomorphism $NS(X)\simeq \RR^d$ for some positive integer $d$ and choose coordinates $x_1,\ldots,x_d$
		on this vector space. Denote by $\Phi_1,\ldots,\Phi_d$ the map $\Phi$ expressed in
		these coordinates and note each $\Phi_j(\be)$ is affine in $\be$. Since we know that 
		${\rm Nef}(X)=\cap_{k=1}^N\{\sum_{j=1}^d a_{kj}x_j\ge 0\}$ with $a_{jk}\in\QQ$ and $N\in\NN$,
		it follows that 
		$\overline{\hbox{\rm AA}(X,D)}=\{\be\in[0,1]^r\,:\, \sum_{j=1}^d a_{kj}\Phi_j(\be)\ge 0,\, k=1,\ldots,N\}$, 
		a set cut out by finitely-many affine equations with rational coefficients in $\be$, proving 
		Theorem \ref{CasciniRThm}.
		
		\eremark

		\begin{spacing}{0}
			
			\def\bi{\bibitem}

		\end{spacing}
		
		\bigskip	
		\bigskip	
		\bigskip	
		\bigskip	
		
		{\sc Imperial College}	
		
		{\tt p.cascini@imperial.ac.uk}
		
		\bigskip	
		
		\bigskip	
		
		{\sc University of Essex}
		
		{\tt jesus.martinez-garcia@essex.ac.uk}
		
		\bigskip	
		\bigskip	
		
		{\sc University of Maryland }
		
		{\tt yanir@alum.mit.edu}

	\end{document}